\documentstyle[11pt,twoside]{article}
\title{\bf A $(p,\nu)$-extension of Srivastava's triple hypergeometric function $H_B$ and its properties}
\author{\sc S.A. Dar$^a$ and R.B. Paris$^b$\\
\\
${}^a\!$ {\em Department of Applied Sciences and Humanities, Faculty of Engineering }\\
{\em and Technology, Jamia Millia Islamia, New Delhi, 110025, India}\\ 
{\em E-Mail: showkatjmi34@gmail.com}\\
${}^b\!$ {\em Division of Computing and Mathematics, Abertay University,}\\
{\em Dundee DD1 1HG, UK}\\
{\em E-Mail: r.paris@abertay.ac.uk}
}
\begin{document}
\newcommand{\bee}{\begin{equation}}
\newcommand{\ee}{\end{equation}}
\newcommand{\br}{\biggr}
\newcommand{\bl}{\biggl}
\newcommand{\g}{\Gamma}
\def\f#1#2{\mbox{${\textstyle \frac{#1}{#2}}$}}
\def\dfrac#1#2{\displaystyle{\frac{#1}{#2}}}
\newcommand{\fr}{\frac{1}{2}}
\newcommand{\fs}{\f{1}{2}}
\date{}
\maketitle
\pagestyle{myheadings}
\markboth{\hfill \it S.A. Dar and R.B. Paris  \hfill}
{\hfill \it A $(p,\nu)$ extension of the $H_B$ function \hfill}
\begin{abstract}
In this paper, we obtain a $(p,\nu)$-extension of Srivastava's triple hypergeometric function $H_B(\cdot)$,
together by using the extended Beta function $B_{p,\nu}(x,y)$ introduced in \cite{AB30}.  We give some of the main properties of this extended function, which include several integral representations involving Exton's hypergeometric function, the Mellin transform, a differential formula, recursion formulas and a bounded inequality. 
\vspace{0.4cm}

\noindent {\bf MSC:} 33C60, 33C65, 33C70, 33B15, 33C05, 33C45, 33C10
\vspace{0.3cm}

\noindent {\bf Keywords:} Srivastava's triple hypergeometric functions; Beta and Gamma functions; Exton's triple hypergeometric function; Bessel function; bounded inequality 
\end{abstract}

\vspace{0.3cm}

\noindent $\,$\hrulefill $\,$

\vspace{0.2cm}

\begin{center}
{\bf 1. \  Introduction and Preliminaries}
\end{center}
\setcounter{section}{1}
\setcounter{equation}{0}
\renewcommand{\theequation}{\arabic{section}.\arabic{equation}}
In the present paper, we employ the following notations:  
\[{\bf N}:=\{1,2,...\},~~{\bf N}_{0}:={\bf N}\cup\{0\},~~{\bf Z}_{0}^{-}:={\bf Z}^{-}\cup\{0\},\]
where the symbols ${\bf N}$ and ${\bf Z}$ denote the set of integer and natural numbers; 
as usual, the symbols ${\bf R}$ and ${\bf C}$  denote the set of real and complex numbers, respectively.

In the available literature, the hypergeometric series and its generalizations appear in various branches of mathematics associated with applications. A large number of triple hypergeometric functions have been introduced and investigated. The work of Srivastava and Karlsson \cite[Chapter 3]{AB37} provides a table of 205 distinct triple hypergeometric functions. Srivastava  introduced the triple hypergeometric functions $H_{A} , H_{B}$ and $H_{C}$ of the second order in \cite{AB34,AB35}. It is known that $H_{C}$ and $H_{B}$ are generalizations of Appell's hypergeometric functions $F_{1}$ and $F_{2}$, while $H_{A}$ is the generalization of both $F_{1}$ and $F_{2}$.

In the present study, we confine our attention to Srivastava's triple hypergeometric function $H_{B}$ given by \cite[p. 43, 1.5(11) to 1.5(13)]{AB37} (see also \cite{AB34} and \cite[p.~68]{AB36})
\[
H_{B}(b_{1},b_{2},b_{3};c_{1},c_{2},c_{3};x,y,z)
:=\sum_{m,n,k=0}^{\infty}\frac{(b_{1})_{m+k}(b_{2})_{m+n}(b_{3})_{n+k}}{(c_{1})_{m}(c_{2})_{n}(c_{3})_{k}}
\frac{x^{m}}{m!}\frac{y^{n}}{n!}\frac{z^{k}}{k!}\]
\bee\label{PV2}
=\sum_{m,n,k=0}^{\infty}\frac{(b_{1}\!+\!b_2)_{2m+n+k}(b_{3})_{n+k}}{(c_{1})_{m}(c_{2})_{n}(c_{3})_{k}}\,\frac{B(b_1\!+\!m\!+\!k,b_2\!+\!m\!+\!n)}{B(b_1,b_2)}\,
\frac{x^{m}}{m!}\frac{y^{n}}{n!}\frac{z^{k}}{k!}.
\ee
Here $(\lambda)_{n}$ denotes the Pochhammer symbol (or the shifted factorial,
 since $(1)_{n}=n!)$  defined by
\[
(\lambda)_n:=\frac{\Gamma(\lambda+n)}{\Gamma(\lambda)}=\left\{\begin{array}{ll}1, & (n=0,\ \lambda\in{\bf C}\backslash\{0\})\\
\lambda(\lambda+1)...(\lambda+n-1),  & (n\in{\bf N},~\lambda\in{\bf C}) 
 \end{array}\right.
\]
and $B(\alpha, \beta)$ denotes the classical Beta function  defined by \cite[(5.12.1)]{AB26}
\begin{equation}\label{ES5}
B(\alpha, \beta)=\left\{\begin{array}{ll} \displaystyle{\int_{0}^{1}t^{\alpha-1}(1-t)^{\beta-1}dt}, & (\Re(\alpha)>0,  \Re(\beta)>0)\\
\\
\displaystyle{\frac{\Gamma(\alpha)\Gamma(\beta)}{\Gamma(\alpha+\beta)}}, & ( (\alpha,\beta)\in{\bf C}\backslash {\bf Z}_{0}^{-}). \end{array}\right.
\end{equation}
The convergence region for Srivastava's triple hypergeometric series $H_{B}(\cdot)$ is given in \cite[ p.243]{ab22}
as $|x|<\alpha$, $|y|<\beta$, $|z|<\gamma$, where $\alpha$, $\beta$, $\gamma$ satisfy the relation $\alpha+\beta+\gamma+2\sqrt{\alpha\beta\gamma}=1$.

A different type of triple hypergeometric function is Exton's  function $ X_{4}(\cdot)$, which is defined by (see \cite{AB19} and \cite[p.~84, Entry (45a)]{AB37})
\begin{equation}\label{PV3}
X_{4}(b_{1},b_{2};c_{1},c_{2},c_{3};x,y,z):=\sum_{m,n,k=0}^{\infty}\frac{(b_{1})_{2m+n+k}(b_{2})_{n+k}}{(c_{1})_{m}(c_{2})_{n}(c_{3})_{k}}
\frac{x^{m}}{m!}\frac{y^{n}}{n!}\frac{z^{k}}{k!}.
\end{equation}
The convergence region for this series is $2\sqrt{|x|}+(\sqrt{|y|}+\sqrt{|z|})^{2}<1$.
We shall also find it convenient to introduce an additional parameter $a$ into $H_B(\cdot)$ in the form
\[H_{B}^{(a)}(b_{1},b_{2},b_{3};c_{1},c_{2},c_{3};x,y,z):=\hspace{6cm}\]
\bee\label{PV2a}
\sum_{m,n,k=0}^{\infty}\frac{(b_{1}\!+\!b_2)_{2m+n+k}(b_{3})_{n+k}}{(c_{1})_{m}(c_{2})_{n}(c_{3})_{k}}\,\frac{B(b_1\!+\!a\!+\!m\!+\!k,b_2\!+\!a\!+\!m\!+\!n)}{B(b_1,b_2)}\,
\frac{x^{m}}{m!}\frac{y^{n}}{n!}\frac{z^{k}}{k!}, 
\ee 
which reduces to (\ref{PV2}) when $a=0$.

In 1997, Chaudhry {\it et al.} \cite[Eq.(1.7)]{AB8} gave a $p$-extension of the Beta function $B(x, y)$ given by \[B(x,y; p)= \int_0^1 t^{x-1}(1-t)^{y-1}\,\exp \left[\frac{-p}{t(1-t)}\right]dt,\qquad (\Re (p)>0)\]
and they proved that this extension has connections with the Macdonald, error and Whittaker functions. Also,  Chaudhry {\it et al.} \cite{AB10} extended the Gaussian hypergeometric series ${}_{2}F_{1}(\cdot)$ and its integral representations. 
Recently, Parmar {\it et al.} \cite{AB30} have given a further extension of the extended Beta function $B(x,y;p)$ by adding one more parameter $\nu$, which we denote and define by
\begin{equation}\label{ES20}
 B_{p,\nu}(x,y)=\sqrt{\frac{2p}{\pi}}\displaystyle{\int_{0}^{1}t^{x-\frac{3}{2}}(1-t)^{y-\frac{3}{2}}}
 K_{\nu+\frac{1}{2}}\left(\frac{p}{t(1-t)}\right)dt,
 \end{equation}
where $\Re(p)>0$, $\nu\geq0$ and $K_{\nu}(z)$ is the modified Bessel function (sometimes known as the Macdonald function) of order $\nu$. When $\nu=0$, (\ref{ES20}) reduces to $B(x,y;p)$, since $K_\frac{1}{2}(z)=\sqrt{\pi/(2z)} e^{-z}$. A different generalization of the Beta function has been given in \cite{AB28}. 

Motivated by some of the above-mentioned extensions of special functions,  many authors
have studied integral representations of $H_{B}(\cdot)$ functions; see  \cite{AB11, AB12, AB13, AB14}.  Our aim in this paper is to introduce a $(p,v)$-extension of Srivastava's triple hypergeometric function $H_{B}(\cdot)$ in (\ref{PV2}), which we denote by $H_{B,p,\nu}(\cdot)$, based on the extended Beta function in (\ref{ES20}), and to systematically investigate some properties of this extended function.  We consider the Mellin transform, a differential formula, recursion formulas and a bounded inequality satisfied by this function.

The plan of this paper as follows. The extended Srivastava hypergeometric function $H_{B,p,\nu}(\cdot)$ is  defined in Section 2 and some integral representations are presented involving the modified Bessel function and Exton's function $X_4$. The main properties of $H_{B,p,\nu}(\cdot)$ namely, its Mellin transform, a differential formula, a bounded inequality and recursion formulas are established in Sections 3-6. Some concluding remarks are made in Section 7.
\vspace{0.6cm}

\begin{center}
{\bf 2. \ The $(p,v)-$extended Srivastava triple hypergeometric function $H_{B,p,\nu}(\cdot)$ }
\end{center}
\setcounter{section}{2}
\setcounter{equation}{0}
\renewcommand{\theequation}{\arabic{section}.\arabic{equation}}
Srivastava introduced the triple hypergeometric function $H_{B}(\cdot)$, together with its integral representations, in \cite{AB34} and \cite{AB36}. Here we consider the following  $(p,\nu)$-extension of this function, which we denote by $H_{B,p,\nu}(\cdot)$, based on the extended Beta
function $B_{p,\nu}(x,y)$ defined in (\ref{ES20}). This is given by
\[H_{B,p,\nu}(b_{1},b_{2},b_{3};c_{1},c_{2},c_{3};x,y,z)\]
\bee\label{PV16} =\sum_{m,n,k=0}^{\infty}\frac{(b_{1}\!+\!b_{2})_{2m+n+k}(b_{3})_{n+k}}{(c_{1})_{m}(c_{2})_{n}(c_{3})_{k}}\frac{B_{p,\nu}(b_{1}\!+\!m\!+\!k,b_{2}\!+\!m\!+\!n)}
{B(b_{1},b_{2})} \frac{x^{m}}{m!}\frac{y^{n}}{n!}\frac{z^{k}}{k!},
\ee
where the parameters $b_{1},b_{2},b_{3}\in{\bf C}$ and  $c_{1},c_{2}, c_{2}\in{\bf C}\backslash {\bf Z}_{0}^{-}$. The region of convergence is $|x|<r$, $|y|<s$, $|z|<t$, where $r+s+t+2\sqrt{rst}=1$.
This definition clearly reduces to the original classical function when $\nu=0$.

Several integral representations for $H_{B,p,\nu}(\cdot)$ involving Exton's triple hypergeometric function in (\ref{PV3}) can be given. We have
\newtheorem{theorem}{Theorem}
\begin{theorem}$\!\!\!.$\ Each of the following integral representations of the extended Srivastava triple hypergeometric function $H_{B,p,\nu}(\cdot)$ holds for $\Re (p)>0$, and $ \min\{\Re(b_{1}),\Re(b_{2})\} > 0$:
\[
H_{B,p,\nu}(b_{1},b_{2},b_{3};c_{1},c_{2},c_{3};x,y,z)
=\frac{\Gamma(b_{1}+b_{2})}{\Gamma(b_{1})\Gamma(b_{2})}\sqrt{\frac{2p}{\pi}}\int_{0}^{1}t^{b_{1}-\frac{3}{2}}(1-t)^{b_{2}-\frac{3}{2}}
K_{\nu+\frac{1}{2}}\left(\frac{p}{t(1-t)}\right)\]
\bee\label{PV17}
\hspace{3cm}\times X_{4}(b_{1}+b_{2},b_{3};c_{1},c_{2},c_{3};xt(1-t),y(1-t),zt)\,dt;
\ee
\[
H_{B,p,\nu}(b_{1},b_{2},b_{3};c_{1},c_{2},c_{3};x,y,z)=\frac{(\beta-\gamma)^{b_{1}-\frac{1}{2}}(\alpha-\gamma)^{b_{2}-\frac{1}{2}}}{(\beta-\alpha)^{b_{1}+b_{2}-2}}
\frac{\Gamma(b_{1}+b_{2})}{\Gamma(b_{1})\Gamma(b_{2})}\sqrt{\frac{2p}{\pi}}\]
\bee\label{PV18}
\times\int_{\alpha}^{\beta}\frac{(\xi-\alpha)^{b_{1}-\frac{3}{2}}(\beta-\xi)^{b_{2}-\frac{3}{2}}}{(\xi-\gamma)^{b_{1}+b_{2}-1}}
K_{\nu+\frac{1}{2}}\left(\frac{p}{\sigma_1 \sigma_2}\right)
X_{4}( b_{1}+b_{2},b_{3};c_{1},c_{2},c_{3}; \sigma_1\sigma_2 x,\sigma_1 y, \sigma_2 z)\,d\xi,
\ee
where 
\[\sigma_{1}=\frac{(\alpha-\gamma)(\beta-\xi)}{(\beta-\alpha)(\xi-\gamma)},\qquad \sigma_{2}=\frac{(\beta-\gamma)(\xi-\alpha)}{(\beta-\alpha)(\xi-\gamma)}\quad (\gamma<\alpha<\beta);\]
\[H_{B,p,\nu}(b_{1},b_{2},b_{3};c_{1},c_{2},c_{3};x,y,z)
=\frac{2\Gamma(b_{1}+b_{2})}{\Gamma(b_{1})\Gamma(b_{2})}\sqrt{\frac{2p}{\pi}}\int_{0}^{\frac{\pi}{2}}(\sin^{2}\xi)^{b_{1}-1}(\cos^{2}\xi)^{b_{2}-1}\]
\bee\label{PV20}
\times
K_{\nu+\frac{1}{2}}\left(\frac{p}{\sigma_1 \sigma_2}\right)X_{4}(b_{1}+b_{2},b_{3};~c_{1},c_{2},c_{3};~\sigma_1 \sigma_2 x,\sigma_1 y, \sigma_2 z)\,d\xi,
\ee
where 
\[\sigma_1=\cos^2 \xi,\qquad \sigma_2=\sin^2 \xi;\]
\[H_{B,p,\nu}(b_{1},b_{2},b_{3};c_{1},c_{2},c_{3};x,y,z)=\frac{2(1+\lambda)^{b_{1}-\frac{1}{2}}\Gamma(b_{1}+b_{2})}{\Gamma(b_{1})\Gamma(b_{2})}
\sqrt{\frac{2p}{\pi}}\int_{0}^{\frac{\pi}{2}}\frac{(\sin^{2}\xi)^{b_{1}-1}(\cos^{2}\xi)^{b_{2}-1}}{(1+\lambda \sin^{2}\xi)^{b_{1}+b_{2}-1}}\]
\bee\label{PV21}
\times K_{\nu+\frac{1}{2}}\left(\frac{p}{\sigma_1 \sigma_2}\right)X_{4}(b_{1}+b_{2},b_{3};c_{1},c_{2},c_{3};\sigma_1 \sigma_2 x,\sigma_1 y, \sigma_2 z)\,d\xi,
\ee
where
\[\sigma_{1}=\frac{\cos^{2}\xi}{1+\lambda \sin^{2}\xi},\qquad\sigma_{2}=\frac{(1+\lambda)\sin^{2}\xi}{1+\lambda \sin^{2}\xi} \quad (\lambda>-1);\]
and
\[H_{B,p,\nu}(b_{1},b_{2},b_{3};c_{1},c_{2},c_{3};x,y,z)
=\frac{2\lambda^{b_{1}-\frac{1}{2}}\Gamma(b_{1}+b_{2})}{\Gamma(b_{1})\Gamma(b_{2})}\]
\[\times\sqrt{\frac{2p}{\pi}}\int_{0}^{\frac{\pi}{2}}\frac{(\sin^{2}\xi)^{b_{1}-1}(\cos^{2}\xi)^{b_{2}-1}}{(\cos^{2}\xi+\lambda \sin^{2}\xi)^{b_{1}+b_{2}-1}}
K_{\nu+\frac{1}{2}}\left(\frac{p}{\sigma_1 \sigma_2}\right)\]
\bee\label{PV23}
\times X_{4}(b_{1}+b_{2},b_{3};c_{1},c_{2},c_{3};\sigma_1 \sigma_2 x,\sigma_1 y, \sigma_2 z)\,d\xi,
\ee
where
\[\sigma_{1}=\frac{\cos^{2}\xi}{\cos^{2}\xi+\lambda \sin^{2}\xi},\qquad \sigma_{2}=\frac{\lambda \sin^{2}\xi}{\cos^{2}\xi+\lambda \sin^{2}\xi} \quad (\lambda>0).
\]
\end{theorem}
\textbf{Proof}:\  The proof of the first integral representation (\ref{PV17}) follows by use of the extended beta function (\ref{ES20}) in (\ref{PV16}), a change the order of integration and summation (with uniform convergence of the integral) and, after simplification, use of Exton's triple hypergeometric function (\ref{PV3}), to obtain the right-hand side of the result (\ref{PV17}). The integral representations (\ref{PV18})-(\ref{PV23}) can be proved directly by using the following transformations
\begin{eqnarray*}
(\ref{PV18}):\qquad t&=&\frac{(\beta-\gamma)(\xi-\alpha)}{(\beta-\alpha)(\xi-\gamma)},\quad \frac{dt}{d\xi}=\frac{(\beta-\gamma)(\alpha-\gamma)}{(\beta-\alpha) (\xi-\gamma)^2},\\
(\ref{PV20}):\qquad t&=& \sin^2 \xi,\quad \frac{dt}{d\xi}=2\sin \xi \cos \xi\\
(\ref{PV21}):\qquad t&=&\frac{(1+\lambda)\sin^{2}\xi}{1+\lambda \sin^{2}\xi},\quad\frac{dt}{d\xi}=\frac{2(1+\lambda)\sin \xi \cos \xi}{(1+\lambda \sin^2 \xi)^2},\\
(\ref{PV23}):\qquad t&=&\frac{\lambda \sin^{2}\xi}{\cos^2{\xi}+\lambda \sin^{2}\xi},\quad \frac{dt}{d\xi}=\frac{2\lambda \sin \xi \cos \xi}{(\cos^2\xi+\lambda \sin^2\xi)^2}\\ 
\end{eqnarray*}
in turn in (\ref{PV17}) to obtain the right-hand side of each result.
\vspace{0.6cm}

\begin{center}
{\bf 3. \  The Mellin transform for $H_{B,p,\nu}(\cdot)$}
\end{center}
\setcounter{section}{3}
\setcounter{equation}{0}
\renewcommand{\theequation}{\arabic{section}.\arabic{equation}}
The Mellin transform of a locally integrable function $f(x)$ on $(0,\infty)$ is given by (see, for example,  \cite[p.193, Section 2.1]{AB125})
\begin{equation}\label{PV27}
\Phi(s)={\cal M}\left\{f(x)\right\}(s)=\int_{0}^{\infty}x^{s-1}f(x)\,dx
\end{equation}
which defines an analytic function in its strip of analyticity $a<\Re (s)<b$.
The inverse Mellin transform of the above function (\ref{PV27}) is defined by
\begin{equation}\label{PV28}
f(x)={\cal M}^{-1}\left\{\Phi(s)\right\}=\frac{1}{2\pi i}\int_{c-i\infty}^{c+i\infty}x^{-s}\Phi(s)\,ds \qquad (a<c<b).
\end{equation}
\begin{theorem}$\!\!\!.$\ The following Mellin transform of the extended Srivastava triple hypergeometric function $ H_{B,p,\nu}(\cdot)$  holds true:
\[
{\cal M}\left\{H_{B,p,\nu}(b_{1},b_{2},b_{3};c_{1},c_{2},c_{3};x,y,z)\right\}(s)=\int_{0}^{\infty}p^{s-1}
H_{B,p,\nu}(b_{1},b_{2},b_{3};c_{1},c_{2},c_{3};x,y,z)\,dp,\]
\begin{equation}\label{PV30}
=\frac{2^{s-1}}{\sqrt{\pi}}\Gamma\left(\frac{s-\nu}{2}\right)\Gamma\left(\frac{s+\nu+1}{2}\right)
H_{B}^{(s)}\left(b_{1},b_{2},b_{3};c_{1},c_{2},c_{3};x,y,z\right),
\end{equation}
where $\Re\,(s)>\nu>0$, $c_{1},c_{2},c_{3}\in{\bf C}\backslash{\bf Z}_{0}^{-}$
and $H_B^{(s)}(\cdot)$ is defined in (\ref{PV2a}).
\end{theorem}
\noindent
\textbf{Proof}: Substituting the extended Srivastava function (\ref{PV16}) into the integral on the left-hand side of (\ref{PV30}) and changing the order of integration (by the uniform convergence of the integral), we obtain
\[{\cal M}\left\{H_{B,p,\nu}(b_{1},b_{2},b_{3};c_{1},c_2,c_3;x,y,z)\right\}(s)\hspace{3cm}\] \[=\sum_{m,n,k=0}^{\infty}\frac{(b_{1}+b_{2})_{2m+n+k}(b_{3})_{n+k}}{(c_{1})_{m}(c_{2})_{n}(c_{3})_{k}~B(b_{1},b_{2})}\frac{x^{m}}{m!}\frac{y^{n}}{n!}\frac{z^{k}}{k!}
  \left\{\int_{0}^{\infty}p^{s-1}B_{p,\nu}(b_{1}+m+k,b_{2}+m+n)\,dp\right\}.
\]
Using the extended Beta function (\ref{ES20}) then shows that
\[{\cal M}\left\{H_{B,p,\nu}(b_{1},b_{2},b_{3};c_{1},c_{2},c_{3};x,y,z)\right\}(s)
=\sqrt{\frac{2}{\pi}}\sum_{m,n,k=0}^{\infty}\frac{(b_{1}+b_{2})_{2m+n+k}(b_{3})_{n+k}}{(c_{1})_{m}(c_{2})_{n}(c_{3})_{k}~B(b_{1},b_{2})}\frac{x^{m}}{m!}\frac{y^{n}}{n!}\frac{z^{k}}{k!}\]
\[\times\int_{0}^{1}t^{b_{1}+m+k-\frac{3}{2}}(1-t)^{b_{2}+m+n-\frac{3}{2}} \left\{\int_{0}^{\infty}p^{s-\frac{1}{2}}K_{\nu+\frac{1}{2}}\left(\frac{p}{t(1-t)}\right)dp\right\}dt.\]
Application of the result \cite[(10.43.19)]{AB26}
\[\int_{0}^{\infty}w^{s-\frac{1}{2}}K_{\alpha+\frac{1}{2}}(w)dw=2^{s-\frac{3}{2}}\Gamma\left(\frac{s-\alpha}{2}\right)
\Gamma\left(\frac{s+\alpha+1}{2}\right)\qquad (|\Re (\alpha)|<\Re (s))\]
followed by the substitution $w=p/(t(1-t))$ produces
\begin{eqnarray*}\label{PV34}
{\cal M}\left\{H_{B,p,\nu}(b_{1},b_{2},b_{3};c_{1},c_{2},c_{3};x,y,z)\right\}(s)
=\frac{2^{s-1}}{\sqrt{\pi}}\Gamma\left(\frac{s-\nu}{2}\right)\Gamma\left(\frac{s+\nu+1}{2}\right)
\nonumber\\ \times\sum_{m,n,k=0}^{\infty}\frac{(b_{1}+b_{2})_{2m+n+k}(b_{3})_{n+k}}{(c_{1})_{m}(c_{2})_{n}(c_{3})_{k}~B(b_{1},b_{2})}\frac{x^{m}}{m!}\frac{y^{n}}{n!}\frac{z^{k}}{k!}
\left\{\int_{0}^{1}t^{b_{1}+m+k+s-1}(1-t)^{b_{2}+m+n+s-1}dt\right\}.
\end{eqnarray*}
Evaluation of the integral in terms of the classical Beta function then finally yields
\begin{eqnarray*}\label{PV35}
\Phi(s)={\cal M}\left\{H_{B,p,\nu}(b_{1},b_{2},b_{3};c_{1},c_{2},c_{3};x,y,z)\right\}(s)
=\frac{2^{s-1}}{\sqrt{\pi}}\Gamma\left(\frac{s-\nu}{2}\right)\Gamma\left(\frac{s+\nu+1}{2}\right)\\  \times\sum_{m,n,k=0}^{\infty}\frac{(b_{1}+b_{2})_{2m+n+k}(b_{3})_{n+k}}{(c_{1})_{m}(c_{2})_{n}(c_{3})_{k}}
\frac{B(b_{1}+m+k+s,b_{2}+m+n+s)}{B(b_{1},b_{2})}\frac{x^{m}}{m!}\frac{y^{n}}{n!}\frac{z^{k}}{k!}.
\end{eqnarray*}
Identifying the above sum as $H_B^{(s)}(b_{1},b_{2},b_{3};c_{1},c_{2},c_{3};x,y,z)$ in (\ref{PV2a}), we obtain the right-hand side of (\ref{PV30}).

\bigskip

\noindent
\textbf{Corollary 1}: \ The following inverse Mellin formula for $H_{B,p,\nu}(\cdot)$ holds:
\[H_{B,p,\nu}(b_{1},b_{2},b_{3};c_{1},c_{2},c_{3};x,y,z)={\cal M}^{-1}\left\{\Phi(s)\right\}\hspace{4cm}\]
\bee\label{PV36}
=\frac{\pi^{-3/2}}{4i}\int_{c-i\infty}^{c+i\infty}\left(\frac{2}{p}\right)^{s}\Gamma\left(\frac{s\!-\!\nu}{2}\right)\Gamma\left(\frac{s\!+\!\nu\!+\!1}{2}\right)
H_{B}^{(s)}\left(b_{1},b_{2},b_{3};c_{1},c_{2},c_{3};x,y,z\right)ds,
\ee
where $c>\nu$.
\vspace{0.6cm}

\begin{center}
{\bf 4. \  A differentiation formula for $H_{B,p,\nu}(\cdot)$}
\end{center}
\setcounter{section}{4}
\setcounter{equation}{0}
\renewcommand{\theequation}{\arabic{section}.\arabic{equation}}
\begin{theorem}$\!\!\!.$\ The following derivative formula for $H_{B,p,\nu}(\cdot)$ holds:
\[\frac{\partial^{M+N+K}}{\partial x^{M}\partial y^{N}\partial z^{K}}~H_{B,p,\nu}\left(b_{1},b_{2},b_{3};c_{1},c_{2},c_{3};x,y,z\right)
=\frac{(b_{1})_{M+K}(b_{2})_{M+N}(b_{3})_{N+K}}{(c_{1})_{M}(c_{2})_{N}(c_{3})_{K}}\]
\bee\label{PV37}
\times H_{B,p,\nu}(b_{1}\!+\!M\!+\!K,b_{2}\!+\!M\!+\!N,b_{3}\!+\!N\!+\!K;c_{1}\!+\!M,c_{2}\!+\!N, c_{3}\!+\!K;x,y,z),
\ee
where $M$, $N$, $K\in{\bf N}_{0}$.
\end{theorem}

\noindent
\textbf{Proof}:\ If we differentiate partially the series for ${\cal H}\equiv H_{B,p,\nu}(b_{1},b_{2},b_{3};c_{1},c_{2},c_{3};x,y,z)$ in (\ref{PV16}) with respect to $x$ we obtain
\[\frac{\partial{\cal H}}{\partial x}
=\sum_{m=1}^\infty \sum_{n,k=0}^\infty \frac{(b_1+b_2)_{2m+n+k} (b_3)_{n+k}}{(c_1)_m(c_2)_n(c_3)_k}\,\frac{B_{p,\nu}(b_1+m+k,b_2+m+n)}{B(b_1, b_2)}\,\frac{x^{m-1}}{(m-1)!}\frac{y^n}{n!} \frac{z^k}{k!}.\]
Making use of the fact that 
\bee\label{eB}
B(b_1,b_2)=\frac{(b_1+b_2)_2}{b_1 b_2}\,B(b_1+1, b_2+1)
\ee
and $(\lambda)_{m+n}=(\lambda)_m(\lambda+m)_n$, we have upon setting $m \to m+1$
\[\frac{\partial{\cal H}}{\partial x}
=\frac{b_1 b_2}{c_1} \sum_{m,n,k=0}^\infty \frac{(b_1+b_2+2)_{2m+n+k} (b_3)_{n+k}}{(c_1+1)_m(c_2)_n(c_3)_k}\,\frac{B_{p,\nu}(b_1\!+\!1\!+\!m\!+\!k,b_2\!+\!1\!+\!m\!+\!n)}{B(b_1+1, b_2+1)}\,\frac{x^{m}}{m!}\frac{y^n}{n!} \frac{z^k}{k!}\]
\bee\label{e41}
=\frac{b_1 b_2}{c_1}\,H_{B,p,\nu}(b_1+1,b_2+1,b_3;c_1+1,c_2,c_3;x,y,z).
\ee
Repeated application of (\ref{e41}) then yields for $M=1, 2, \ldots $
\[\frac{\partial^M}{\partial x^M}{\cal H}=
\frac{(b_1)_M (b_2)_M}{(c_1)_M}\,H_{B,p,\nu}(b_{1}+M,b_{2}+M,b_{3};c_{1}+M,c_{2},c_{3};x,y,z).\]

A similar reasoning shows that
\[\frac{\partial^{M+1}}{\partial x^M \partial y}{\cal H}
=\frac{(b_1)_M (b_2)_M}{(c_1)_M} \sum_{m,k=0}^\infty \sum_{n=1}^\infty \frac{(b_1\!+\!b_2\!+\!2M)_{2m+n+k} (b_3)_{n+k}}{(c_1\!+\!M)_m (c_2)_n(c_3)_k}\hspace{3cm}\]
\[\hspace{3cm}\times\frac{B_{p, \nu}(b_1\!+\!M\!+\!m\!+\!k,b_2\!+\!M\!+\!n\!+\!k)}{B(b_1\!+\!M, b_2\!+\!M)}\,\frac{x^m}{m!} \frac{y^{n-1}}{(n-1)!}\frac{z^k}{k!}\]
\bee\label{e42}
=\frac{(b_1)_{M} (b_2)_{M+1} b_3}{(c_1)_{M}c_2}\,H_{B,p,\nu}(b_1+M,b_2+M+1,b_3+1;c_1+M,c_2+1,c_3;x,y,z)
\ee
upon putting $n\to n+1$ and using the property of the Beta function in (\ref{ES5}).
Repeated differentiation of (\ref{e42}) $N$ times with respect to $y$ then produces 
\[\frac{\partial^{M+N}}{\partial x^M \partial y^N} {\cal H}=\frac{(b_1)_M (b_2)_{M+N} (b_3)_N}{(c_1)_M (c_2)_N}\,
H_{B,p,\nu}(b_1\!+\!M,b_2\!+\!M\!+\!N,b_3\!+\!N;c_1\!+\!M,c_2\!+\!N,c_3;x,y,z).\]

Application of the same procedure to deal with differentiation with respect to $z$ then yields the result stated in (\ref{PV37}).

\vspace{0.6cm}

\begin{center}
{\bf 5. \  An upper bound for $H_{B,p,\nu}(\cdot)$}
\end{center}
\setcounter{section}{5}
\setcounter{equation}{0}
\renewcommand{\theequation}{\arabic{section}.\arabic{equation}}
\begin{theorem}$\!\!\!.$\ Let the parameters $b_j$, $c_j$ $(1\leq j\leq 3)$ be positive and the variables $x$, $y$, $z$ be complex. Then the following bounded inequality for $H_{B,p,\nu}(\cdot)$ holds:
\[|H_{B,p,\nu}(b_{1},b_{2},b_{3};c_{1},c_2,c_3;x,y,z)|\hspace{8cm}\]
\bee\label{e51}
<\frac{2^\nu |p|^{\nu+1}}{\sqrt{\pi} (\Re (p))^{2\nu+1}}\,\Gamma(\nu+\fs)\,H_B^{(\nu)}(b_1,b_2,b_3;c_1,c_2,c_3;|x|,|y|,|z|),
\ee
where $\Re(p)>0$, $\nu>0$ and $H_B^{(\nu)}(\cdot)$ is defined in (\ref{PV2a})
\end{theorem}

The integral representation of the extension $H_{B,p,\nu}(\cdot)$ in (\ref{PV17}) is associated with the modified Bessel function of the second kind, for which we have the following expression \cite[(10.32.8)]{AB26}
\[ K_{\nu+\frac{1}{2}}(z)=\frac{\sqrt{\pi}\left(\frac{1}{2}z\right)^{\nu+\frac{1}{2}}}{\Gamma(\nu+1)}\int_{1}^{\infty}e^{-zt}(t^{2}-1)^{\nu}dt,\qquad(\nu>-1,\ \Re (z)>0).
\]
In our problem we have $\nu>0$, $\Re(z)>0$. Further, we let $x=\Re(z)$, so that
\[| K_{\nu+\frac{1}{2}}(z)|\leq\frac{\sqrt{\pi}\left(\frac{1}{2}|z|\right)^{\nu+\frac{1}{2}}}{\Gamma(\nu+1)}\left|\int_{1}^{\infty}e^{-zt}(t^{2}-1)^{\nu}dt\right|<\frac{\sqrt{\pi}\left(\frac{1}{2}|z|\right)^{\nu+\frac{1}{2}}}{\Gamma(\nu+1)} \int_0^1 t^{2\nu}e^{-xt} dt \]
\begin{equation}\label{ES87}
=\frac{\sqrt{\pi}\left(\frac{1}{2}|z|\right)^{\nu+\frac{1}{2}}}{\Gamma(\nu+1)}\frac{\Gamma(2\nu+1,x)}{x^{2\nu+1}},
\end{equation}
where $\Gamma(a,z)$ is the upper incomplete gamma function \cite[(8.2.2)]{AB26}. Although this bound is numerically found to be quite sharp when $z $ is real, it involves the incomplete gamma function which would make the integral for $F_{1,p,\nu}(b_{1},b_{2},b_{3};c_{1};x,y)$ difficult to bound. We can simplify (\ref{ES87}) by making use of the simple inequality $\Gamma(2\nu+1,x)<\Gamma(2\nu+1)$ to find
\begin{equation}\label{ES88}
| K_{\nu+\frac{1}{2}}(z)|<\frac{\sqrt{\pi}\left(\frac{1}{2}|z|\right)^{\nu+\frac{1}{2}}}{\Gamma(\nu+1)}\frac{\Gamma(2\nu+1)}{x^{2\nu+1}}
=\frac{1}{2}\left(\frac{2|z|}{x^{2}}\right)^{\nu+\frac{1}{2}}\Gamma(\nu+\fs),
\end{equation}
upon use of the duplication formula for the gamma function. The bound (\ref{ES88}) is less sharp than (\ref{ES87}) but has the advantage of being easier to handle in the integral for $H_{B,p,\nu}(\cdot)$. 
\\
\\
\noindent
\textbf{Proof}: Setting $z=p/(t(1-t))$, where $t\in(0,1)$ and $\Re(p)>0$, in (\ref{ES88}) we obtain
\[
\left| K_{\nu+\frac{1}{2}}\left(\frac{p}{t(1-t)}\right)\right|<\frac{1}{2}\left(\frac{2|p|t(1-t)}{(\Re (p))^{2}}\right)^{\nu+\frac{1}{2}}\Gamma(\nu+\fs).
\]
For ease of presentation we shall assume that the parameters $b_j$, $c_j>0$ ($1\leq j\leq 3$). Then, from (\ref{PV17}),
\[|H_{B,p,\nu}(b_1, b_2, b_3;c_1,c_2,c_3;x,y,z)|
\leq \frac{\sqrt{2|p|/\pi}}{B(b_1,b_2)}\int_0^1\bigg|t^{b_1-\frac{3}{2}}(1-t)^{b_2-\frac{3}{2}} \,K_{\nu+\frac{1}{2}}\left(\frac{p}{t(1-t)}\right)\]
\[\times X_4(b_1+b_2,b_3;c_1,c_2,c_3;xt(1-t),y(1-t),zt) \bigg|\,dt\]
\[<\frac{2^\nu |p|^{\nu+1}}{\sqrt{\pi} (\Re (p))^{2\nu+1}}\,\frac{\Gamma(\nu+\fs)}{B(b_1,b_2)} \sum_{m,n,k=0}^\infty \frac{(b_1+b_2)_{2m+n+k} (b_3)_{n+k}}{(c_1)_m(c_2)_n(c_3)_k}\,\frac{|x|^m}{m!}\frac{|y|^n}{n!}\frac{|z|^k}{k!}\]
\[\hspace{5cm}\times\int_0^1t^{b_1+\nu+m+k-1}(1-t)^{b_2+\nu+m+n-1} dt\]
\bee\label{e52}
<\frac{2^\nu |p|^{\nu+1}\Gamma(\nu+\fs)}{\sqrt{\pi} (\Re (p))^{2\nu+1}}\,\sum_{m,n,k=0}^\infty 
\frac{(b_1+b_2)_{2m+n+k}(b_3)_{n+k}}{(c_1)_m(c_2)_n(c_3)_k}\hspace{4.5cm}\]
\[\hspace{3cm}\times\frac{B(b_1\!+\!\nu\!+\!m\!+\!k,b_2\!+\!\nu\!+\!m\!+\!n)}{B(b_1,b_2)}\,\frac{|x|^m}{m!}\frac{|y|^n}{n!}\frac{|z|^k}{k!}
\ee
which is the result stated in (\ref{e51}).
 
 \vspace{0.6cm}

\begin{center}
{\bf 6. \  Recursion formulas for $H_{B,p,\nu}(\cdot)$}
\end{center}
\setcounter{section}{6}
\setcounter{equation}{0}
\renewcommand{\theequation}{\arabic{section}.\arabic{equation}}
In this section, we obtain two recursion formulas for the extended Srivastava function $H_{B,p,\nu}(\cdot)$. The first formula gives a recursion with respect to the numerator parameter $b_3$, and the second a recursion with respect to any one of the denominator parameters $c_j$ ($1\leq j\leq 3$).
\begin{theorem}$\!\!\!.$\ The following recursion for $H_{B,p,\nu}(\cdot)$ with respect to the numerator parameter $b_3$ holds:
\[H_{B,p,\nu}(b_1,b_2,b_3\!+\!1;c_1,c_2,c_3;x,y,z)=H_{B,p,\nu}(b_1,b_2,b_3;c_1,c_2,c_3;x,y,z)\]
\bee\label{e60}
+\frac{yb_2}{c_2} H_{B,p,\nu}(b_1,b_2\!+\!1,b_3\!+\!1;c_1, c_2+1,c_3; x,y,z)+\frac{zb_1}{c_3} H_{B,p,\nu}(b_1\!+\!1,b_2,b_3\!+\!1;c_1, c_2,c_3\!+\!1; x,y,z).
\ee
\end{theorem}
\noindent
{\bf Proof.}\ \ From (\ref{PV16}) and the result $ (b_{3}+1)_{n+k}=(b_{3})_{n+k}(1+n/b_{3}+k/b_{3})$, we obtain
\[
H_{B,p,\nu}(b_{1},b_{2},b_{3}+1;c_{1},c_{2},c_{3};x,y,z)\hspace{7cm}\] \[=\sum_{m,n,k=0}^{\infty}\frac{(b_{1}+b_{2})_{2m+n+k}(b_{3}+1)_{n+k}}{(c_{1})_{m}(c_{2})_{n}(c_{3})_{k}}\frac{B_{p,\nu}(b_{1}+m+k,b_{2}+m+n)}
{B(b_{1},b_{2})} \frac{x^{m}}{m!}\frac{y^{n}}{n!}\frac{z^{k}}{k!}\]
\[
=H_{B,p,\nu}(b_{1},b_{2},b_{3};c_{1},c_{2},c_{3};x,y,z)\hspace{7cm}\]
\[+\frac{y}{b_{3}}\sum_{m=0}^{\infty}\sum_{n=1}^{\infty}\sum_{k=0}^{\infty}\frac{(b_{1}+b_{2})_{2m+n+k}(b_{3})_{n+k}}{(c_{1})_{m}(c_{2})_{n}(c_{3})_{k}}\frac{B_{p,\nu}(b_{1}+m+k,b_{2}+m+n)}
{B(b_{1},b_{2})} \frac{x^{m}}{m!}\frac{y^{n-1}}{(n-1)!}\frac{z^{k}}{k!}\]
\bee\label{e61}
+\frac{z}{b_{3}}\sum_{m=0}^{\infty}\sum_{n=0}^{\infty}\sum_{k=1}^{\infty}\frac{(b_{1}+b_{2})_{2m+n+k}(b_{3})_{n+k}}{(c_{1})_{m}(c_{2})_{n}(c_{3})_{k}}\frac{B_{p,\nu}(b_{1}+m+k,b_{2}+m+n)}
{B(b_{1},b_{2})} \frac{x^{m}}{m!}\frac{y^{n}}{n!}\frac{z^{k-1}}{(k-1)!}.
\ee

Consider the first sum in (\ref{e61}) which we denote by $S$. Put $n\to n+1$ and use the identity $(a)_{n+1}=a (a+1)_n$ to find
\[S=\frac{y}{b_3}\sum_{m,n,k=0}^\infty  \frac{(b_1\!+\!b_2)_{2m+n+1+k} (b_3)_{n+1+k}}{(c_1)_m (c_2)_{n+1} (c_3)_k}\,
\frac{B_{p,\nu}(b_1\!+\!m\!+\!k,b_2\!+\!1\!+\!m\!+\!n)}{B(b_1, b_2)}\,\frac{x^m}{m!} \frac{y^{n}}{n!} \frac{z^k}{k!}\] 
\[=\frac{y(b_1\!+\!b_2)}{c_2}\sum_{m,n,k=0}^\infty  \frac{(b_1\!+\!b_2\!+\!1)_{2m+n+k} (b_3+1)_{n+k}}{(c_1)_m (c_2+1)_{n} (c_3)_k}\,
\frac{B_{p,\nu}(b_1\!+\!m\!+\!k,b_2\!+\!1\!+\!m\!+\!n)}{B(b_1, b_2)}\,\frac{x^m}{m!} \frac{y^{n}}{n!} \frac{z^k}{k!}.\]  
Using the fact that $$B(b_1, b_2)=\frac{b_1+b_2}{b_2} B(b_1,b_2+1),$$ we then obtain
\[S=\frac{yb_2}{c_2}\sum_{m,n,k=0}^\infty  \frac{(b_1\!+\!b_2\!+\!1)_{2m+n+k} (b_3\!+\!1)_{n+k}}{(c_1)_m (c_2\!+\!1)_{n} (c_3)_k}\,
\frac{B_{p,\nu}(b_1\!+\!m\!+\!k,b_2\!+\!1\!+\!m\!+\!n)}{B(b_1, b_2+1)}\,\frac{x^m}{m!} \frac{y^{n}}{n!} \frac{z^k}{k!}\] 
\bee\label{e62}
=\frac{yb_2}{c_2} H_{B,p,\nu}(b_1,b_2+1,b_3+1;c_1, c_2+1,c_3; x,y,z).
\ee

Proceeding in a similar manner for the second series in (\ref{e61}) with $k\to k+1$, we find that this sum can be expressed as
\bee\label{e63}
\frac{zb_1}{c_3} H_{B,p,\nu}(b_1+1,b_2,b_3+1;c_1, c_2,c_3+1; x,y,z).
\ee
Combination of (\ref{e62}) and (\ref{e63}) with (\ref{e61}) then produces the result stated in (\ref{e60}).
\bigskip

\noindent
\textbf{Corollary 2}: \ From (\ref{e60}) the following recursion holds
\[H_{B,p,\nu}(b_1,b_2,b_3+N;c_1,c_2,c_3;x,y,z)=H_{B,p,\nu}(b_1,b_2,b_3;c_1,c_2,c_3;x,y,z)\]
\[+\frac{yb_2}{c_2}\sum_{\ell=1}^N H_{B,p,\nu}(b_1,b_2+1,b_3+\ell;c_1,c_2+1,c_3;x,y,z)\]
\bee\label{e63a}
+\frac{zb_1}{c_3}\sum_{\ell=1}^NH_{B,p,\nu}(b_1+1,b_2,b_3+\ell;c_1,c_2,c_3+1;x,y,z)
\ee
for positive integer $N$.

\begin{theorem}$\!\!\!.$\ The following 3-term recursion for $H_{B,p,\nu}(\cdot)$ with respect to the denominator parameter $c_1$ holds:
\[H_{B,p,\nu}(b_1,b_2,b_3;c_1,c_2,c_3;x,y,z)=\]
\bee\label{e64}
H_{B,p,\nu}(b_1,b_2,b_3;c_1+1,c_2,c_3;x,y,z)+\frac{x b_1 b_2}{c_1(c_1+1)}
H_{B,p,\nu}(b_1+1,b_2+1,b_3;c_1+2,c_2,c_3;x,y,z).
\ee
Permutation of the $c_j$ enables analogous recursions in the denominator parameters $c_2$ and $c_3$ to be obtained.
\end{theorem}
\bigskip

\noindent
{\bf Proof.}\ \ 
Consider the case when $c_1$ is reduced by 1, namely
\[H\equiv H_{B,p,\nu}(b_1,b_2,b_3;c_1-1,c_2,c_3;x,y,z)\]
and use $(c_1-1)_m=(c_1)_m/\{1+\frac{m}{c_1-1}\}$.
Then 
\[H=\sum_{m,n,k=0}^\infty\frac{(b_1\!+\!b_2)_{2m+n+k} (b_3)_{n+k}}{(c_1\!-\!1)_m (c_2)_{n} (c_3)_k}\,
\frac{B_{p,\nu}(b_1\!+\!m\!+\!k,b_2\!+\!m\!+\!n)}{B(b_1, b_2)}\,\frac{x^m}{m!} \frac{y^{n}}{n!} \frac{z^k}{k!}\]
\[=\sum_{m,n,k=0}^\infty\frac{(b_1\!+\!b_2)_{2m+n+k} (b_3)_{n+k}}{(c_1)_m (c_2)_{n} (c_3)_k}\,
\frac{B_{p,\nu}(b_1\!+\!m\!+\!k,b_2\!+\!m\!+\!n)}{B(b_1, b_2)}\,\bl(1+\frac{m}{c_1-1}\br)\frac{x^m}{m!} \frac{y^{n}}{n!} \frac{z^k}{k!}\]
\[=H_{B,p,\nu}(b_1,b_2,b_3;c_1,c_2,c_3;x,y,z)\hspace{8cm}\]
\[+\frac{x}{c_1-1}\sum_{m=1}^\infty\sum_{n,k=0}^\infty\frac{(b_1\!+\!b_2)_{2m+n+k} (b_3)_{n+k}}{(c_1)_m (c_2)_{n} (c_3)_k}\,
\frac{B_{p,\nu}(b_1\!+\!m\!+\!k,b_2\!+\!m\!+\!n)}{B(b_1, b_2)}\,\frac{x^{m-1}}{(m-1)!} \frac{y^{n}}{n!} \frac{z^k}{k!}.\]
Putting $m\to m+1$ in the above sum, we obtain
\[\frac{x}{c_1-1}\sum_{m,n,k=0}^\infty\frac{(b_1+b_2)_{2m+2+n+k} (b_3)_{n+k}}{(c_1)_{m+1} (c_2)_{n} (c_3)_k}\,
\frac{B_{p,\nu}(b_1\!+\!1\!+\!m\!+\!k,b_2\!+\!1\!+\!m\!+\!n)}{B(b_1, b_2)}\,\frac{x^{m}}{m!} \frac{y^{n}}{n!} \frac{z^k}{k!}\]
\[=\frac{x(b_1\!+\!b_2)_2}{c_1(c_1-1)}\sum_{m,n,k=0}^\infty\frac{(b_1\!+\!b_2\!+\!2)_{2m+n+k} (b_3)_{n+k}}{(c_1+1)_{m} (c_2)_{n} (c_3)_k}\,
\frac{B_{p,\nu}(b_1\!+\!1\!+\!m\!+\!k,b_2\!+\!1\!+\!m\!+\!n)}{B(b_1, b_2)}\,\frac{x^{m}}{m!} \frac{y^{n}}{n!} \frac{z^k}{k!}.
\]
Using (\ref{eB}), we find that this last sum becomes
\[\frac{x b_1 b_2}{c_1(c_1-1)}\sum_{m,n,k=0}^\infty\frac{(b_1\!+\!b_2\!+\!2)_{2m+n+k} (b_3)_{n+k}}{(c_1+1)_{m} (c_2)_{n} (c_3)_k}\,
\frac{B_{p,\nu}(b_1\!+\!1\!+\!m\!+\!k,b_2\!+\!1\!+\!m\!+\!n)}{B(b_1+1, b_2+1)}\,\frac{x^{m}}{m!} \frac{y^{n}}{n!} \frac{z^k}{k!}\]
\[=\frac{x b_1 b_2}{c_1(c_1-1)}\,H_{B,p,\nu}(b_1+1,b_2+1,b_3;c_1+1,c_2,c_3;x,y,z).\]

This then yields the recurrence relation (in $c_1$) given by
\[H_{B,p,\nu}(b_1,b_2,b_3;c_1-1,c_2,c_3;x,y,z)=\]
\[H_{B,p,\nu}(b_1,b_2,b_3;c_1,c_2,c_3;x,y,z)+\frac{x b_1 b_2}{c_1(c_1-1)}
H_{B,p,\nu}(b_1+1,b_2+1,b_3;c_1+1,c_2,c_3;x,y,z).\]
Replacement of $c_1$ by $c_1+1$ then yields the result stated in (\ref{e64}).

\vspace{0.6cm}

\begin{center}
{\bf 7. \  Concluding remarks}
\end{center}
\setcounter{section}{7}
\setcounter{equation}{0}
\renewcommand{\theequation}{\arabic{section}.\arabic{equation}}
In this paper, we have introduced the $(p,\nu)$-extended Srivastava triple hypergeometric function given by $H_{B,p,\nu}(\cdot)$  in (\ref{PV16}). We have given some integral representations of this function that involve the modified Bessel function of the second kind and Exton's triple hypergeometric function $X_4$. We have also established some properties of the function $H_{B,p,\nu}(\cdot)$, namely the Mellin transform, a differential formula, a bounded inequality  and some recursion relations. 


\vspace{0.6cm}

\begin{thebibliography}{99}
\footnotesize{
\bibitem{AB6} Bozer, M. and \"{O}zarslan, M.A.,  Notes on generalized gamma, beta and hypergeometric functions. J. Comput. Anal. Appl. \textbf{15(7)} (2013) 1194-1201.

\bibitem{AB8} Chaudhry, M.A., Qadir, A., Rafique, M. and  Zubair, S.M., Extension of Euler's beta function. J. Comput. Appl. Math., \textbf{78(1)} (1997) 19-32.
 
\bibitem{AB10} Chaudhry, M.A., Qadir, A., Srivastava, H.M. and Paris, R.B., Extended hypergeometric and confluent hypergeometric functions. Appl. Math. Comput., \textbf{159(2)} (2004) 589-602.
 
\bibitem{AB9}  \c{C}etinkaya, A., Ya\u{g}basana, M.I. and Kiymaz, I.O., The extended Srivastava's triple hypergeometric functions and their integral representations. J. Nonlinear Sci. Appl. (2016) 1-1.

\bibitem{AB11} Choi, J., Hasanov, A.,  Srivastava, H.M.  and Turaev, M., Integral representations for Srivastava's triple hypergeometric functions. Taiwanese J. Math. \textbf{15(6)} (2011) 2751-2762.

\bibitem{AB12} Choi, J., Hasanov, A. and Turaev, M., Integral representations for Srivastava's triple hypergeometric functions $H_{A}$. Honam Mathematical J. \textbf{34(1)} (2012) 113-124.

\bibitem{AB13} Choi, J., Hasanov, A. and  Turaev, M., Integral representations for Srivastava's triple hypergeometric functions $H_{B}$. J. Korean Soc. Math. Educ. Ser. B: Pure Appl. Math. \textbf{19(2)} (2012) 137-145.

\bibitem{AB14} Choi, J., Hasanov, A. and Turaev, M., Integral representations for Srivastava's triple hypergeometric functions $H_{C}$. Honam Mathematical J. \textbf{34(4)} (2012) 473-482.

\bibitem{AB15} Choi, J., Rathie, A.K. and Parmar, R.K., Extension of extended Beta, hypergeometric function and confluent hypergeometric function. Honam Mathematical J. \textbf{36(2)} (2014) 357-385.

\bibitem{AB16}  Choi, J., Parmar, R.K.  Po\'{g}any, T.K., Mathieu-type series built by $(p, q)$ extended Gaussian hypergeometric function. arXiv:1604.05077v1 [math.CA] (2016) 1-9.

\bibitem{AB17} Erd\'{e}lyi, A., Magnus, W., Oberhettinger, F. and Tricomi, F.G., \textit{Higher Transcendental Functions}. Vol. 1. McGraw-Hill, New York, Toronto and London, 1953.


\bibitem{AB19} Exton, H., Hypergeometric functions of three variables. J. Indian Acad. Math., {\bf 4} (1982) 113-119.


\bibitem{ab22} Karlsson, P.W., Regions of convergences for hypergeometric series in three variables. Math. Scand, \textbf{48} (1974) 241-248.


\bibitem{AB24}  Luo, M.J., Milovanovic, G.V. and  Agarwal, P., Some results on the extended beta and extended hypergeometric functions. Applied Math. and Computation  \textbf{248} (2014) 631-651.



\bibitem{AB125} Oberhettinger, F., \textit{Tables of Mellin Transforms.} Springer-Verlag, Berlin, Heidelberg, New York, 1974.


\bibitem{AB26}Olver, F.W.J., Lozier, D.W., Boisvert, R.F. and Clark C.W. (eds.), \textit{NIST Handbook of Mathematical Functions}. Cambridge University Press, Cambridge, 2010.


\bibitem{AB28}  \"{O}zarslan, M.A.,  Some remarks on extended hypergeometric, extended confluent hypergeometric and extended Appell's functions. J. Comput. Anal. Appl. \textbf{14(6)} (2012) 1148-1153.

\bibitem{AB29}  \"{O}zergin, E.,   \"{O}zarslan, M.A. and Altin, A., Extension of gamma, beta and hypergeometric functions. J. Comput. Appl. Math. \textbf{235} (2011) 4601-4610.

\bibitem{AB30} Parmar, R.K., Chopra, P.  and Paris, R.B., On an Extension of extended Beta and hypergeometric Functions. arXiv:1502.06200 [math.CA] \textbf{22} (2015). [to appear in J. Classical Anal.]


\bibitem{AB34} Srivastava, H.M., Hypergeometric functions of three variables. Ganita \textbf{15} (1964) 97-108.

\bibitem{AB35} Srivastava, H.M.,  Some integrals representing triple hypergeometric functions. Rend. Circ. Mat. Palermo (Ser. 2) \textbf{16} (1967) 99-115.

\bibitem{AB36} Srivastava, H.M. and Manocha, H.L., \textit{A Treatise on Generating Functions} Halsted Press (Ellis Horwood Limited, Chichester, U.K.) John Wiley and Sons, New York, Chichester, Brisbane and Toronto, 1984.

\bibitem{AB37} Srivastava, H.M. and Karlsson, P.W., \textit{Multiple Gaussian Hypergeometric Series}. Halsted Press (Ellis Horwood Limited, Chichester), John Wiley and Sons, New York, Chichester, Brisbane and Toronto, 1985.

}
\end{thebibliography}
\end{document}